%% file: main_arXiv_v2.tex
\pdfoutput=1

\documentclass[3p,authoryear, times]{elsarticle}
\journal{Statistics and Probability Letters}

\include{packages}

\usepackage{appendix}

\include{macros}

\begin{document}
	
\begin{frontmatter}
\title{Uniform concentration bounds for frequencies of rare events}

\author[1]{Stéphane Lhaut}
\ead{stephane.lhaut@uclouvain.be}

\author[2]{Anne Sabourin}
\ead{anne.sabourin@telecom-paris.fr}

\author[1]{Johan Segers}
\ead{johan.segers@uclouvain.be}

\address[1]{Institut de Statistique, Biostatistique et Sciences Actuarielles, UCLouvain, Voie du Roman Pays 20, 1348
	Louvain-la-Neuve, Belgium}
\address[2]{LTCI, Télécom Paris, Institut polytechnique de Paris, 19 place Marguerite Perey, 91120 Palaiseau, France}

\begin{abstract}
New Vapnik–Chervonenkis type concentration inequalities are derived for the empirical distribution of an independent
random sample. Focus is on the maximal deviation over classes of Borel sets within a low probability region. The
constants are explicit, enabling numerical comparisons.
\end{abstract}
	
\end{frontmatter}
	
\section{Introduction}
\label{sec: introduction}
	
\input{introduction}

\section{Inequalities for rare events}
\label{sec: extreme}

\input{extreme}

\section{Working directly with the maximal deviation}
\label{sec: maxdev}

\input{maxdev}

\section{Working with the expected maximal deviation}
\label{sec: expectation}

\input{expectation}

\section{Comparison of the bounds}
\label{sec: comparison}

\input{comparison}

\bibliography{references.bib}

\newpage
\appendix
\gdef\thesection{\Alph{section}} 
\makeatletter
\renewcommand\@seccntformat[1]{Appendix \csname the#1\endcsname.\hspace{0.5em}}
\makeatother

\appendixpage

\section{An improvement of the Vapnik--Chervonenkis inequality}
\label{sec: an improvement}

\input{anImprovement}

\section{Vapnik--Chervonenkis inequality for relative deviations}
\label{sec: vc relatives}

\input{VCRelatives}

\section{Proof of Theorem~\ref{thm: VC ineq rare sym before}}
\label{sec: proof sym before}

\input{proofSymBefore}

\end{document}

%% file: packages.tex
\usepackage[english]{babel}
\usepackage[T1]{fontenc}
\usepackage[utf8]{inputenc}

\usepackage{amssymb}
\usepackage{amsmath}
\usepackage{amsthm}
\usepackage{amscd}
\usepackage{centernot} 
\usepackage{mathtools}
\usepackage{IEEEtrantools} 

\usepackage{graphicx}
\usepackage{array}
\usepackage{listings}
\usepackage{color}
\usepackage{xcolor}
\usepackage{comment} 
\usepackage{enumitem}
\usepackage{dsfont}
\usepackage{float}
\usepackage{csquotes}
\usepackage{caption}
\usepackage{textcomp}

\usepackage{natbib}
\biboptions{authoryear}

\usepackage[colorlinks,menucolor=blue,linkcolor=blue, citecolor=blue, urlcolor=blue]{hyperref} 

\allowdisplaybreaks

%% file: macros.tex
\newcommand{\R}{\mathbb{R}}
\newcommand{\N}{\mathbb{N}}

\newcommand{\ra}{\rightarrow}

\newcommand{\e}{\mathrm{e}}

\newcommand{\lp}{\left(}
\newcommand{\rp}{\right)}
\newcommand{\lc}{\left[}
\newcommand{\rc}{\right]}
\newcommand{\lacc}{\left\{}
\newcommand{\racc}{\right\}}

\DeclareMathOperator{\E}{E}
\DeclareMathOperator{\Var}{Var}

\renewcommand{\P}{\operatorname{P}}
\newcommand{\1}{\mathds{1}}

\DeclarePairedDelimiter\floor{\lfloor}{\rfloor}
\newcommand{\A}{\mathcal{A}}
\renewcommand{\S}{\mathbb{S}}

\DeclareMathOperator{\Bin}{\mathit{Bin}}

\numberwithin{equation}{section}

\theoremstyle{plain}
\newtheorem{theoreme}{Theorem}[section]
\newtheorem{proposition}[theoreme]	 {Proposition}	
\newtheorem{corollary}	[theoreme]	 {Corollary}	
\newtheorem{lemma}	    [theoreme]	 {Lemma}
	
\theoremstyle{definition}

\theoremstyle{remark}
\newtheorem{remark}	    [theoreme]	 {Remark}

%% file: introduction.tex
Let $X_1,\ldots,X_n$ be an independent random sample with common distribution $\mu$ on $\R^d$. We always let $\mu_n$ denote the associated \emph{empirical measure}, 
i.e., $\mu_n := n^{-1} \sum_{i=1}^n \delta_{X_i}$.
We already know, from \citet[Theorem~2]{Vapnik1971}, the \emph{Vapnik and Chervonenkis (VC) inequality} which states that 
\begin{equation}
\label{eq: VC concentration}
	\P \lp \sup_{A \in \A} |\mu_n(A) - \mu(A)| \geq t \rp \leq 4 \S_\A(2n)\e^{-nt^2/8},
\end{equation}
where $\A$ is a non-empty class of Borel sets on $\R^d$ and $\S_\A(n) := \max_{x_1,\ldots,x_n \in \R^d} \big| \big\{ \{x_1,\ldots,x_n\} \cap A : A \in \A \big\} \big|$ is the shattering coefficient of the class $\A$ -- for background on VC theory, we refer to~\citet{Lugosi2002}.
Setting the upper bound equal to $\delta \in (0,1)$ and solving for $t$ yields that, with probability $1-\delta$,
\begin{equation}
\label{eq: VC ineq}
	\sup_{A \in \A} |\mu_n(A) - \mu(A)| \leq 2 \sqrt{\frac{2}{n} \big[ \log(4/\delta) + \log \S_\A(2n) \big]}.
\end{equation}

In Appendix~\ref{sec: an improvement}, we show an improvement of this inequality that will be used in Section~\ref{sec: extreme}: with probability at least $1-\delta$,
\[
	\sup_{A \in \A} |\mu_n(A) - \mu(A)| \leq \sqrt{\frac{1}{2n}} + \sqrt{\frac{2}{n} \big[ \log(4/\delta) + \log \S_\A(2n) \big]}.
\]

Our main goal is to develop bounds with \emph{explicit constants} when working with \emph{low probability regions} (Section~\ref{sec: extreme}), i.e., it will be assumed that there exists some Borel set $\mathbb{A} \subseteq \R^d$ which contains every set $A \in \A$ and has a “small” probability 
\begin{equation}
\label{eq: rare events}
	p:=\mu(\mathbb{A}).
\end{equation}
For example, consider an empirical risk minimization problem where we use the empirical risk to classify some data. In that situation, $p$ could be $1\%$ if we are looking to understand how the prediction behaves when the input variables take values above the 99th percentile of the distribution of $\|X\|$, where $X$ has the same distribution as the data in the training set -- we refer to \citet[Remark~5]{Goix2015}, \citet{Jalalzai2018} and \citet[Section~4.2]{Clemencon2021} for classification in extreme regions.  

The bounds should get smaller since all the encountered events have a small probability; however, the estimation is more difficult since we have less data in those regions. In this rare events framework, a first improvement of the VC inequality \eqref{eq: VC ineq}, based on~\citet[Theorem~2.1]{Anthony1993} and~\citet[Theorem~1.11]{Lugosi2002}, is given by
\begin{equation}
\label{eq: VC ineq Anthony}
	\sup_{A \in \A} |\mu_n(A) - \mu(A)| \leq 2 \sqrt{\frac{2p}{n} \big[ \log(12/\delta) + \log \S_\A(2n) \big]}. 
\end{equation}
The details of the statement and the proof are deferred to Appendix~\ref{sec: vc relatives}.
Even though the factor $p$ in the square root will improve the bound, this result is still not fully satisfactory. The \emph{effective sample size}, i.e., the expected number of data points in the low probability region, is $np$ rather than $n$. Therefore, the shattering coefficient involved in the bound seems too large. 

This work is heavily motivated by~\citet{Goix2015} who already introduced a VC type inequality adapted to rare events. However, unlike the classical VC bound, the constants appearing in their result are \emph{not} explicit. A similar remark can be made about~\citet{Gine2006} where concentration inequalities are also introduced for normalized empirical processes, in a very general setting. In Section~\ref{sec: maxdev} and Section~\ref{sec: expectation}, we discuss different methods that lead to new, explicit, inequalities. Essentially, there are two possibilities: either we directly apply the standard tools of VC theory on the maximal deviation itself (Section~\ref{sec: maxdev}), or we use a variation of the classical \emph{McDiarmid bounded differences inequality}~\citep[Theorem~3.8]{McDiarmid1998} and then we deal with the expectation $\E \lc \sup_{A \in \A} |\mu_n(A) - \mu(A)| \rc$ (Section~\ref{sec: expectation}). After deriving the different bounds, we compare them in Section~\ref{sec: comparison}.

Our aim is not to be exhaustive in the bounds that we derive. Other tools can lead to other inequalities. As an example, in~\citet[Theorem~A.1]{Clemencon2021}, we provide a bound on the expectation of the supremum using chaining techniques, which leads to a better asymptotic rate, but which is inaccurate for realistic sample sizes due to the large constants inherent to those methods. Therefore, such results are not presented here. A detailed comparison can be found in~\citet{Lhaut2021}.

\begin{remark}[Pointwise measurability]
\label{rem: pm}
To ensure measurability of the supremum appearing in~\eqref{eq: VC concentration}, a common hypothesis is to assume that the class $\A$ is \emph{pointwise measurable}, as suggested by~\citet[Example~2.3.4]{VVV1996}, i.e., that there exists a subclass $\A_0 \subseteq \A$, at most countable, such that for every $A \in \A$, there exists a sequence $(A_n)_{n \in \N} \subseteq \A_0$ such that
\[
	\lim_{n \ra \infty} \1_{A_n}(x) = \1_A(x), \quad \text{for every } x \in \R^d.
\]
It is easily showed that, under this assumption, $\sup_{A \in \A} |\mu_n(A) - \mu(A)| = \sup_{A \in \A_0} |\mu_n(A) - \mu(A)|$ and we may replace $\A$ by $\A_0$ everywhere.
\end{remark}

%% file: extreme.tex
The main tool to develop adapted bounds is expressed in the following lemma, which already appears in~\citet[Page~17]{Goix2015} and the idea of which is likely to be found in other places in the literature as well, see, e.g., \citet[Equation~(14.6)]{Novak2011} for a result in the same spirit. The proof is elementary and can be found in~\citet[Chapter~3]{Lhaut2021}.

\begin{lemma}[Conditioning trick]
\label{lem: cond trick}
Let $X_1,\ldots,X_n$ be an independent random sample with common distribution $\mu$ on $\R^d$ and $\mu_n$ be the associated empirical measure. Let $\A$ be a non-empty class of Borel sets on $\R^d$. Let $\mathbb{A}$ and $p$ be as in~\eqref{eq: rare events}. Let $Y_1,\ldots,Y_n$ be an independent random sample with common distribution 
\[
	\mu_{\mathbb{A}}(\cdot) = \mu(\cdot \, | \, \mathbb{A}) = \frac{\mu(\cdot \cap \mathbb{A})}{\mu(\mathbb{A})} = \frac{\mu(\cdot \cap \mathbb{A})}{p},
\]
and independent of $X_1,\ldots,X_n$. If, for every $k=1,\ldots,n$ we denote $\mu_k^Y$ the empirical measure associated to $Y_1,\ldots,Y_k$ and if $K \sim \Bin(n,p)$ denotes the number of data points $X_i$ in $\mathbb{A}$, then the following equality in distribution holds 
\[
	\lc (\mu_n(A))_{A \in \A} \, \Big| \, K = k \rc \stackrel{d}{=} \lp \frac{k}{n} \mu_k^Y(A) \rp_{A \in \A},
\]
in the sense of equality of finite-dimensional distributions. In particular, under the pointwise measurability assumption (Remark~\ref{rem: pm}), this equality still holds when considering the supremum over the class $\A$. 
\end{lemma}

This result provides us with a convenient way to adapt classical VC inequalities to our setting: we start by conditioning $\mu_n$ on $K$, then we apply concentration inequalities on the empirical measure $\mu_K^Y$ which takes account of the rare nature of the encountered events and finally we integrate out on $K$. In this last step, we will make use of the \emph{Bernstein's inequality} for binomial random variables \cite[Theorem~1.5]{Lugosi2002} which states that for every $t>0$,
\begin{equation}
\label{eq: bernstein ineq}
	\P \big( |K-np| \geq t \big) \leq 2 \exp \lp - \frac{t^2}{2(np(1-p)+t/3)} \rp.
\end{equation}

When we will be working with the expectation, we shall also make use of another version of the VC inequality~\eqref{eq: VC ineq} for the expected maximal deviation:
\begin{equation}
\label{eq: VC ineq mean}
	\E \lc \sup_{A \in \A} |\mu_n(A) - \mu(A)| \rc \leq \sqrt{\frac{2 \log(2\S_\A(2n))}{n}},
\end{equation}
the proof of which can be found, for example, in~\citet[Theorem~1.9]{Lugosi2002}.

In the bounds that we propose, we will need to define the shattering coefficient $\S_\A(x)$ for a \emph{real} parameter $x>0$. It is simply understood that $\S_\A(x) :=  \S_\A(\floor{x})$, where $\floor{x}$ denotes the integer part of $x$.

%% file: maxdev.tex
Using the improvement of the VC inequality that we show in Appendix~\ref{sec: an improvement}, we derive a first inequality. 

\begin{theoreme}
\label{thm: VC ineq rare sym after}
Let $X_1,\ldots,X_n$ be an independent random sample with common distribution $\mu$ on $\R^d$ and let $\mu_n$ be the associated empirical measure. Let $\A$ be a non-empty class of Borel sets on $\R^d$. Let $\mathbb{A}$ and $p$ be as in~\eqref{eq: rare events}. If $np \geq 4\log(4/\delta)$, where $\delta \in (0,1)$, then we have, with probability $1-\delta$,
\[
	\sup_{A \in \A} \left| \mu_n(A) - \mu(A) \right|
	\leq \frac{2}{3n} \log(4/\delta) + \sqrt{\frac{p}{n}} \lp \sqrt{2\log(4/\delta)} + 2 \sqrt{\log(8/\delta) + \log \S_\A(4np)} + 1 \rp.  
\]
\end{theoreme}

\begin{proof}[Proof of Theorem~\ref{thm: VC ineq rare sym after}]
	Let $t>0$. If $K=\sum_{i=1}^n \1\{X_i \in \mathbb{A}\} \sim \Bin(n,p)$, then by Lemma~\ref{lem: cond trick},
	\begin{align*}
		&\P \lp \sup_{A \in \A} \frac{|\mu_n(A) - \mu(A)|}{p} \geq t \rp 
		= \E \lc \P \lp \sup_{A \in \A} \frac{|\mu_n(A) - \mu(A)|}{p} \geq t \, \Big| \, K \rp \rc 
		= \E \lc \P \lp \sup_{A \in \A} \frac{\left| \tfrac{K}{n} \mu_K^Y(A) - \mu(A) \right|}{p} \geq t \, \Big| \, K \rp \rc \\
		&\qquad = \E \lc \P \lp \sup_{A \in \A} \left| \tfrac{K}{np} \mu_K^Y(A) - \mu_{\mathbb{A}}(A) \right| \geq t \, \Big| \, K \rp \rc 
		\leq \E \lc \P \lp \frac{K}{np} \sup_{A \in \A} \left| \mu_K^Y(A) - \mu_\mathbb{A}(A) \right| + \left| \frac{K}{np} - 1 \right| \geq t \, \Big| \, K \rp \rc,
	\end{align*}
	where the last inequality follows from the triangle inequality and the fact that $\mu_{\mathbb{A}}(A) \leq 1$ for any $A \in \A$.
	Let $t=u+s$ for $u,s >0$. Then, the latter quantity is bounded by
	\begin{equation}
		\label{eq: mass separation}
		\P \lp \left| \frac{K}{np} - 1 \right| \geq s \rp + \E \lc \P \lp \frac{K}{np} \sup_{A \in \A} \left| \mu_K^Y(A) - \mu_\mathbb{A}(A) \right| \geq u \, \Big| \, K \rp \1 \lacc \left| \frac{K}{np} - 1 \right| \leq s \racc \rc.
	\end{equation}
	Let $\delta \in (0,1)$. We will pick $u=u(\delta)$ and $s=s(\delta)$ such that each term in \eqref{eq: mass separation} is bounded by $\delta/2$. 
	
	We deal with the first term of \eqref{eq: mass separation} directly using the Bernstein inequality \eqref{eq: bernstein ineq},
	\[
		\P \lp \left| \frac{K}{np} - 1 \right| \geq s \rp  \leq \P \big( \left| K - np \right| \geq nps \big) \leq 2 \exp \lp - \frac{n^2 p^2 s^2}{2(np(1-p)+nps/3)} \rp \leq 2 \exp \lp - \frac{nps^2}{2(1+s/3)} \rp. 
	\]
	The positive root $s_+$ associated with the quadratic equation obtained by equaling this last term to $\delta/2$ satisfies
	\[
	s_+ \leq \frac{2}{3np} \log(4/\delta) + \sqrt{\frac{2}{np} \log(4/\delta)} =: s(\delta).
	\] 
	
	To deal with the second term of \eqref{eq: mass separation}, we use the improved VC inequality that we develop in Appendix~\ref{sec: an improvement} in the form \eqref{eq: improved VC concentration form}, 
	\begin{align*}
		\P \lp \frac{K}{np} \sup_{A \in \A} \left| \mu_K^Y(A) - \mu_\mathbb{A}(A) \right| \geq u \, \Big| \, K \rp
		\leq 4 \S_\A(2K) \exp \lp - \frac{(np)^2u^2}{2K} \Big( 1 - \sqrt{\frac{2K}{(np)^2u^2}} \Big) - \frac{1}{4} \rp.     
	\end{align*}
By monotonicity in $K$ on the region of interest \footnote{The map $x \in (0,+\infty) \mapsto -x(1-\sqrt{1/x}) = -x + \sqrt{x}$ is decreasing on $(1/4, +\infty)$. Our focus is on that region of the real line since, in our situation, $x=\frac{(np)^2u^2}{2K}$ with $u=u(\delta) \geq 1/\sqrt{np}$ and $K \leq np(1+s(\delta)) \leq 2np$.}, we get
	\[
		\E \lc \P \lp \frac{K}{np} \sup_{A \in \A} \left| \mu_K^Y(A) - \mu_\mathbb{A}(A) \right| \geq u \, \Big| \, K \rp \1 \lacc \left| \frac{K}{np} - 1 \right| \leq s \racc \rc \leq  4 \S_\A(2np(1+s)) \exp \lp - \frac{npu^2}{2(1+s)} \Big( 1 - \sqrt{\frac{2(1+s)}{npu^2}} \Big) - \frac{1}{4} \rp.
	\]
	If $np \geq 4\log(4/\delta)$, we have $s(\delta) \leq \frac{1}{6} + \frac{1}{\sqrt{2}} \leq 1$. Hence, using $s=1$ in the latter expression, we find that it equals $\delta/2$ if 
	\[
	u = u(\delta) := \frac{1}{\sqrt{np}} + 2 \sqrt{\frac{1}{np} \big[ \log(8/\delta) + \log \S_\A(4np) \big]}.
	\]
	
	Regrouping the two terms, we obtain that with probability at least $1-\delta$,
	\[
		\sup_{A \in \A} \frac{\left| \mu_n(A) - \mu(A) \right|}{p} \leq \frac{2}{3np} \log(4/\delta) + \sqrt{\frac{1}{np}} \lp \sqrt{2\log(4/\delta)} + 2 \sqrt{\log(8/\delta) + \log \S_\A(4np)} + 1 \rp. 
	\]
	Multiplying both sides by $p$, we get the result.
\end{proof}

Another possibility consists of symmetrizing the process, based on our improved version of the classical argument (Lemma~\ref{lem: symmetrization}), before using the conditioning trick. It leads to a simpler bound. The proof is deferred to Appendix~\ref{sec: proof sym before}.

\begin{theoreme}
	\label{thm: VC ineq rare sym before}
	Let $X_1,\ldots,X_n$ be an independent random sample with common distribution $\mu$ on $\R^d$ and let $\mu_n$ be the associated empirical measure. Let $\A$ be a non-empty class of Borel sets on $\R^d$. Let $\mathbb{A}$ and $p$ be as in~\eqref{eq: rare events}. If $np \geq 2\log(8/\delta)$, where $\delta \in (0,1)$, then we have, with probability $1-\delta$,
	\[
	\sup_{A \in \A} \left| \mu_n(A) - \mu(A) \right| \leq \sqrt{\frac{2p}{n}} \lp 2 \sqrt{\log(8/\delta) + \log \S_\A(8np)} + 1 \rp.  
	\]
\end{theoreme}

%% file: expectation.tex
The variation of the McDiarmid bounded differences inequality that we use is recalled in~\citet[Proposition~11]{Goix2015}. When combined with~\citet[Lemma~12]{Goix2015}, it leads to a convenient Bernstein type concentration inequality for the maximal deviation that is particularly adapted when working with rare events, i.e., for every $t>0$,
\[
\P \lp \sup_{A \in \A} |\mu_n(A) - \mu(A)| - \E \lc \sup_{A \in \A} |\mu_n(A) - \mu(A)| \rc \geq t \rp \leq \exp \lp - \frac{nt^2}{4p+2t/3} \rp.
\]
Setting the upper bound equal to $\delta \in (0,1)$ and solving for $t>0$ proves the following result.

\begin{proposition}[Concentration of the maximal deviation]
	\label{prop: concentration for md}
	Let $X_1,\ldots,X_n$ be an independent random sample with common distribution $\mu$ on $\R^d$ and let $\mu_n$ be the associated empirical measure. Let $\A$ be a non-empty class of Borel sets on $\R^d$. Let $\mathbb{A}$ and $p$ be as in~\eqref{eq: rare events}. Then, for every $\delta \in (0,1)$, with probability at least $1-\delta$, we have
	\[
	\sup_{A \in \A} \left| \mu_n(A) - \mu(A) \right| \leq \frac{2}{3n} \log(1/\delta) + 2 \sqrt{\frac{p}{n}\log(1/\delta)} + \E \lc \sup_{A \in \A} |\mu_n(A) - \mu(A)| \rc.
	\]
\end{proposition}

Using this inequality, we may directly work with the expectation $\E \lc \sup_{A \in \A} |\mu_n(A) - \mu(A)| \rc$ to obtain a bound with high probability for $\sup_{A \in \A} \left| \mu_n(A) - \mu(A) \right|$. 

We start by applying the conditioning trick (Lemma~\ref{lem: cond trick}) on the expectation. 

\begin{lemma}
	\label{lem: cond trick for the expectation}
	Let $X_1,\ldots,X_n$ be an independent random sample with common distribution $\mu$ on $\R^d$ and let $\mu_n$ be the associated empirical measure. Let $\A$ be a non-empty class of Borel sets on $\R^d$. Let $\mathbb{A}$ and $p$ be as in~\eqref{eq: rare events}. Then,
	\[
	\E \lc \sup_{A \in \A} |\mu_n(A) - \mu(A)| \rc \leq \E \lc \frac{K}{n} \sup_{A \in \A} |\mu_K^Y(A) - \mu_{\mathbb{A}}(A)| \rc + \sqrt{\frac{p}{n}},
	\]
	where $\mu_{\mathbb{A}}, \mu_K^Y$ and $K$ are as in Lemma~\ref{lem: cond trick}.
\end{lemma}

\begin{proof}[Proof~of~Lemma~\ref{lem: cond trick for the expectation}]
	By Lemma~\ref{lem: cond trick} and a computation similar to the one in the beginning of the proof of Theorem~\ref{thm: VC ineq rare sym after}, we have
	\begin{align*}
		&\E \lc \sup_{A \in \A} \left| \mu_n(A)-\mu(A) \right| \rc 
		= \E \lc \E \lc \sup_{A \in \A} \left| \mu_n(A)-\mu(A) \right| \, \Big| \, K \rc \rc \\
		&\qquad = \E \lc \sup_{A \in \A} \left| \frac{K}{n} \big( \mu_K^Y(A) - \mu_{\mathbb{A}}(A) \big) + \lp \frac{K}{n} - p \rp \mu_{\mathbb{A}}(A) \right| \rc 
		\leq  \E \lc \frac{K}{n} \sup_{A \in \A} \left| \mu_K^Y(A)-\mu_{\mathbb{A}}(A) \right| \rc + \E \lc \left| \frac{K}{n} - p \right| \rc.
	\end{align*}
	The second term is easily bounded using Jensen's inequality and the fact that $K \sim \Bin(n,p)$:
	\[
	\lp \E \lc \left| \frac{K}{n} - p \right| \rc \rp^2 \leq \E \lc \lp \frac{K}{n} - p \rp^2 \rc = \Var \lc \frac{K}{n} - p \rc = \frac{\Var[K]}{n^2} = \frac{np(1-p)}{n^2} \leq \frac{p}{n}. \qedhere
	\]
\end{proof}

To deal with the remaining expectation $\E \lc \frac{K}{n} \sup_{A \in \A} |\mu_K^Y(A) - \mu_{\mathbb{A}}(A)| \rc$, we will make use of~\eqref{eq: VC ineq mean}. Furthermore, we will also assume that the VC dimension $V_\A := \sup \{n \in \N : \S_\A(n) = 2^n\}$ of the class $\A$ is finite to be able to apply the famous \emph{Sauer's Lemma}~\cite[Lemma~1]{Bousquet2004} which states that for every $n \in \N$,
\begin{equation}
\label{eq: sauer}
	\S_\A(n) \leq (n+1)^{V_\A}.
\end{equation}
The purpose of this assumption is that, whenever $V_\A < +\infty$, we will be able to use Jensen's inequality which is sharper than the monotonicity arguments underlying every proof in Section~\ref{sec: maxdev}.

\begin{proposition}
\label{prop: mean + Jensen}
Let $X_1,\ldots,X_n$ be an independent random sample with common distribution $\mu$ on $\R^d$ and let $\mu_n$ be the associated empirical measure. Let $\A$ be a non-empty class of Borel sets on $\R^d$. Let $\mathbb{A}$ and $p$ be as in~\eqref{eq: rare events}. Then, if $V_\A < +\infty$,
	\[
	\E \lc \frac{K}{n} \sup_{A \in \A} |\mu_K^Y(A) - \mu_{\mathbb{A}}(A)| \rc \leq \sqrt{\frac{2p}{n}\big[ \log2 + V_\A \log(2np+1) \big]}.
	\]
\end{proposition}

\begin{proof}[Proof~of~Proposition~\ref{prop: mean + Jensen}]
	By the VC inequality for the expectation~\eqref{eq: VC ineq mean} combined with Sauer's Lemma~\eqref{eq: sauer} and Jensen's inequality, we have
	\begin{align*}
		&\E \lc \frac{K}{n} \sup_{A \in \A} \left| \mu_K^Y(A)-\mu_{\mathbb{A}}(A) \right| \rc
		 = \E \lc \E \lc \frac{K}{n} \sup_{A \in \A} \left| \mu_K^Y(A)-\mu_{\mathbb{A}}(A) \right| \, 					\Big| \, K \rc \rc \\
		&\quad \leq \E \lc \frac{K}{n} \sqrt{\frac{2 \log \lp 2 \S_\A(2K) \rp}{K}} \rc 
		 \leq \frac{1}{n} \E \lc \sqrt{2K \big[ \log(2) +  V_\A \log(2K+1) \big]} \rc \\
		&\quad \leq \frac{1}{n} \sqrt{2 \E[K] \big[ \log(2) +  V_\A \log(2 \E[K] +1) \big]} 
		 = \sqrt{\frac{2p}{n} \big[ \log(2) +  V_\A \log(2np +1) \big]},   
	\end{align*}
	where we made use of the concavity of the map $K \mapsto \sqrt{2K  \big[ \log(2) +  V_\A \log(2K+1) \big]}$ (the derivative is clearly decreasing) and the fact that $K \sim \Bin(n,p)$. 
\end{proof}

Combining all the ingredients of this section, we finally get a bound on the maximal deviation.
\begin{corollary}
\label{cor: mean + Jensen}
Let $X_1,\ldots,X_n$ be an independent random sample with common distribution $\mu$ on $\R^d$ and let $\mu_n$ be the associated empirical measure. Let $\A$ be a non-empty class of Borel sets on $\R^d$. Let $\mathbb{A}$ and $p$ be as in~\eqref{eq: rare events}. Then, if $V_\A < +\infty$, we have, for any $\delta \in (0,1)$, with probability at least $1-\delta$,
\[
	\sup_{A \in \A} \left| \mu_n(A) - \mu(A) \right| \leq \frac{2}{3n} \log(1/\delta) + \sqrt{\frac{2p}			{n}} \lp \sqrt{2\log(1/\delta)} + \sqrt{\log 2 + V_\A \log(2np+1)} + \frac{\sqrt{2}}{2} \rp.
\]
\end{corollary}

%% file: comparison.tex
We obtained various bounds on the maximal deviation $\sup_{A \in \A} |\mu_n(A) - \mu(A)|$ which hold with probability $1-\delta \in (0,1)$. Before comparing them, we  summarize the main results of our paper and we recall the bound derived in Theorem~\ref{thm: relative deviations}:
\begin{itemize}
	\item Theorem~\ref{thm: VC ineq rare sym after} (working directly on the maximal deviation and  symmetrizing the process after the application of the conditioning trick): whenever $np \geq 4\log(4/\delta)$,
	\begin{equation}
	\label{eq: first bound}
		\sup_{A \in \A} \left| \mu_n(A) - \mu(A) \right| 
		\leq \frac{2}{3n} \log(4/\delta) + \sqrt{\frac{2p}{n}} \lp \sqrt{\log(4/\delta)} + \sqrt{2[\log(8/\delta) + \log \S_\A(4np)]} + 1 \rp.
	\end{equation}
	\item Theorem~\ref{thm: VC ineq rare sym before} (working directly on the maximal deviation and symmetrizing the process before the application of the conditioning trick): 	whenever $np \geq 2\log(8/\delta)$,
	\begin{equation}
	\label{eq: second bound}
		\sup_{A \in \A} \left| \mu_n(A) - \mu(A) \right| \leq \sqrt{\frac{2p}{n}} \lp 2 \sqrt{\log(8/\delta) + \log \S_\A(8np)} + 1 \rp.
	\end{equation}
	\item Corollary~\ref{cor: mean + Jensen} (working on the expected maximal deviation and using Jensen's inequality): 	whenever $V_\A < +\infty$,
    \begin{equation}
    \label{eq: third bound}
		\sup_{A \in \A} \left| \mu_n(A) - \mu(A) \right| \leq \frac{2}{3n} \log(1/\delta) + \sqrt{\frac{2p}{n}} \lp \sqrt{2\log(1/\delta)} + \sqrt{\log 2 + V_\A \log(2np+1)} + \frac{\sqrt{2}}{2} \rp.
	\end{equation}
	\item Theorem~\ref{thm: relative deviations} (based on \citet[Theorem~2.1]{Anthony1993}):  whenever $np \geq 8\log(3/\delta)/3$,
	\begin{equation}
	\label{eq: vc relative}
		\sup_{A \in \A} |\mu_n(A) - \mu(A)| \leq 2 \sqrt{\frac{2p}{n} \big[ \log(12/\delta) + \log 					\S_\A(2n) 	\big]}.
	\end{equation}
\end{itemize}

We observe that the bound~\eqref{eq: third bound} seems to beat~\eqref{eq: first bound} in any case. Indeed, both inequalities have the same structure: one term decaying like $n^{-1}$ and one term decreasing like $n^{-1/2}$; however, those terms seem always lower in the concentration bound obtained through the use of McDiarmid's inequality. Even though it could be that the shattering coefficient of $\A$ grows slower than the rate induced by the use of Sauer's Lemma~\eqref{eq: sauer} in~\eqref{eq: third bound}, the fact that~\eqref{eq: sauer} enables us to use Jensen's inequality heavily compensates the loss. 

The comparison between~\eqref{eq: second bound} and~\eqref{eq: third bound} is a bit harder since the structure is not the same anymore. Nevertheless, the coefficient of the leading term in~\eqref{eq: second bound} seems larger than the one in~\eqref{eq: third bound}. Hence, even though the inequality obtained in Section~\ref{sec: expectation} could be worse for small samples due to its additional term, we think that for practical sample sizes, it remains the best bound.

\begin{figure}[h]
	\centering
	\includegraphics[scale=0.5]{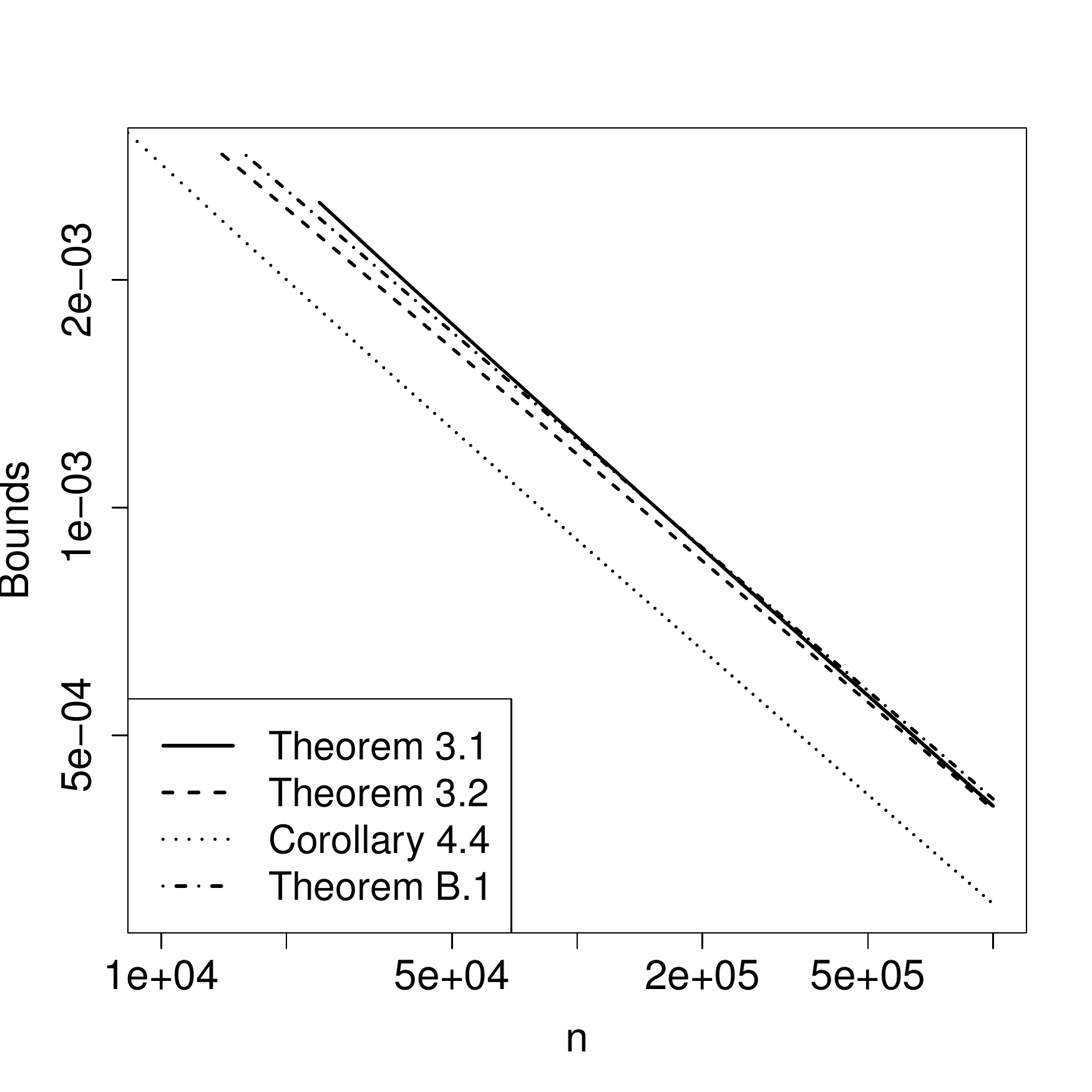}
	\caption{Comparison of the bounds for the empirical cumulative distribution function on sufficiently large samples (to satisfy the hypothesis of the different theorems). Axes are on logarithmic scale.}
	\label{fig: comp}
\end{figure}

In Figure~\ref{fig: comp}, we provide a graphical comparison of our results with the more classical VC inequality for relative deviations~\eqref{eq: vc relative} on the simplest, one-dimensional, VC class that we may think off: $\A = \{ (-\infty,t] : t \in \R,\, t \leq Q(10^{-3})\}$, where $Q(p) := \inf \{ x \in \R : F(x) \geq p\}$ is the left-continuous inverse of $F(x) := \mu((-\infty,x]))$, the cumulative distribution function associated to $\mu$. In this situation, we easily verify that $\S_A(n) = n+1$ and that the main assumption of our paper~\eqref{eq: rare events} is satisfied with $p=10^{-3}$. Note also that $\A$ is pointwise measurable (Remark~\ref{rem: pm}). We used $\delta = 10^{-2}$. 

One clearly observes better performance of our bound~\eqref{eq: third bound} in this situation, which is valid for every sample size. The bounds coming from Section~\ref{sec: maxdev} are closer to the older bound; however, it seems that for samples with a size larger than $10^6$ data points, they perform a bit better. More importantly, those bounds take into account the reduced effective sample size $np$ in their measure of size complexity $\S_\A$, which is not the case with~\eqref{eq: vc relative}.

%% file: anImprovement.tex
We show an improvement of the VC inequality~\cite[Theorem~2.17]{Lhaut2021}.

\begin{theoreme}
\label{thm: improved VC}
Let $X_1,\ldots,X_n$ be an independent sample with common distribution $\mu$ on $\R^d$ and let $\mu_n$ be the associated empirical distribution. Consider a non-empty class $\A$ of Borel sets of $\R^d$. Then, for every $\delta \in (0,1)$, with probability $1-\delta$, we have 
\[
	\sup_{A \in \A} |\mu_n(A) - \mu(A)| \leq \sqrt{\frac{1}{2n}} + \sqrt{\frac{2}{n} \big[ \log(4/\delta) + \log \S_\A(2n) \big]}.
\]
\end{theoreme}

The order of the bound is the same as in the original VC bound~\eqref{eq: VC ineq}. However, the constant in front of the square root is almost halved. The proof relies on a variation of a classical \emph{symmetrization} argument \cite[Lemma~2~p.~185]{Bousquet2004}. 

\begin{lemma}[Symmetrization]
\label{lem: symmetrization}
Let $X_1,\ldots,X_n,X_1',\ldots,X_n'$ be an independent sample with common distribution $\mu$ on $\R^d$. Let $\mu_n$ denote the empirical distribution associated to $X_1,\ldots,X_n$ and let $\mu_n'$ denote the empirical distribution associated to $X_1',\ldots,X_n'$. Consider a non-empty class $\A$ of Borel sets of $\R^d$. If $t > 0$ and $a \in (0,1)$ are such that $4(1-a)^2nt^2 \geq 2$, then 
\[
	\left.
	\begin{array}{ll}
		& \P \big( \sup_{A \in \mathcal{A}} (\mu_n(A) - \mu(A)) \geq t \big)\\[+10pt]
		& \P \big( \sup_{A \in \mathcal{A}} (\mu(A) - \mu_n(A)) \geq t \big)
	\end{array}
	\right \} \leq  2 \P \lp  \sup_{A \in \mathcal{A}} (\mu_n(A) - \mu_n'(A)) \geq at \rp.
\]
\end{lemma}

The sample $X_1',\ldots,X_n'$ is often referred to as a “ghost sample”. The choice $a=1/2$ gives back the initial symmetrization. However, as we will see, this more flexible result will give a way to substantially improve the classical bound. 

\begin{proof}[Proof~of~Lemma~\ref{lem: symmetrization}]
Let $A \in \A$ and $t>0$. We first observe that
\[
	\1{\{\mu_n(A) - \mu(A) \geq t\}} \1{\{\mu_n'(A) - \mu(A) \leq (1-a)t\}} \leq \1{\{\mu_n(A) - \mu_n'(A) \geq at\}}. 
\]
Integrating both sides with respect to the ghost sample, we get
\[
    \1{\{\mu_n(A) - \mu(A) \geq t\}} \P' \big( \mu_n'(A) - \mu(A) \leq (1-a)t \big) \leq \P' \big( \mu_n(A) - \mu_n'(A) \geq at \big),
\]
where we used the notation $\P'$ to emphasize that the integration was realized only with respect to $X_1',\ldots,X_n'$. By Bienaymé--Tchebychev inequality and hypothesis, we deduce
\begin{align*}
	\P' \big( \mu_n'(A) - \mu(A) \leq (1-a)t \big) 
	&= 1 - \P' \big( \mu_n'(A) - \mu(A) \geq (1-a)t \big) \\
	&\geq 1 - \frac{\Var[\mu_n'(A)]}{(1-a)^2t^2} = 1 - \frac{\Var \lc \1{\{X_1' \in A\}}\rc}					{n(1-a)^2t^2} \geq 1 - \frac{1}{4n(1-a)^2t^2} \geq 1/2,
\end{align*}
where we used the fact that the variance of $\1{\{X_1' \in A\}}$ is bounded by $1/4$ (since $\max_{p \in [0,1]} \{p(1-p)\}$ = 1/4). We obtain
\begin{equation*}
	\1{\{\mu_n(A) - \mu(A) \geq t\}} \leq 2 \P' \big( \mu_n(A) - \mu_n'(A) \geq at \big).
\end{equation*}
Since this relation holds for every $A \in \A$, we finally get that
\begin{align*}
	\1{\{\sup_{A \in \mathcal{A}} \mu_n(A) - \mu(A) \geq t\}} &= \sup_{A \in \mathcal{A}}  \1{\{ \mu_n(A) - \mu(A) \geq t\}} 
	\leq 2 \sup_{A \in \mathcal{A}} \P' \big( \mu_n(A) - \mu_n'(A) \geq at \big) \\
	&\leq 2 \P' \big( \sup_{A \in \mathcal{A}} (\mu_n(A) - \mu_n'(A)) \geq at \big). 
\end{align*}
Taking the expectation with respect to the sample $X_1,\ldots,X_n$, it follows that
\begin{equation*}
	\P \big( \sup_{A \in \mathcal{A}} (\mu_n(A) - \mu(A)) \geq t \big) \leq 2 \P \big( \sup_{A \in \mathcal{A}} (\mu_n(A) - \mu_n'(A)) \geq at \big).
\end{equation*}

The other inequality has a similar proof.
\end{proof}

\begin{proof}[Proof~of~Theorem~\ref{thm: improved VC}]
Let $nt^2 \geq 1/2$. Following the same arguments as the ones used to prove the classical VC inequality~\cite[Theorem~2.14]{Lhaut2021} and using our symmetrization, we show that for any $a\in (0,1)$, whenever $4(1-a)^2nt^2 \geq 2$, we have
\[
	\P \lp \sup_{A \in \A} |\mu_n(A) - \mu(A)| \geq t \rp \leq 4 \S_\A(2n)\e^{-a^2nt^2/2}.
\]
Choosing the largest $a \in (0,1)$ such that this relation holds, i.e. $a=a(n,t)=1-\sqrt{1/(2nt^2)}$, we get
\begin{equation}
\label{eq: improved VC concentration form}
	\P \lp \sup_{A \in \A} \left| \mu_n(A) - \mu(A) \right| \geq t \rp \leq 4 \S_\A(2n) \exp{\lacc -\frac{nt^2}{2} \Big( 1 - \sqrt{\frac{2}{nt^2}} \Big) - \frac{1}{4} \racc}. 
\end{equation}
Setting the upper bound equal to $\delta \in (0,1)$ and solving for $t$ leads to a quadratic equation for $t$ whose positive solution corresponds to the bound proposed in the theorem.

If $nt^2 < 1/2$, the bound~\eqref{eq: improved VC concentration form} is trivial.
\end{proof}

%% file: VCRelatives.tex
\begin{theoreme}[VC inequality for relative deviations]
\label{thm: relative deviations}
Let $X_1,\ldots,X_n$ be an independent random sample with common distribution $\mu$ on $\R^d$ and $\mu_n$ be the associated empirical measure. Let $\A$ be a non-empty class of Borel sets on $\R^d$. Let $\mathbb{A}$ and $p$ be as in~\eqref{eq: rare events}. If $np \geq (8/3)\log(3/\delta)$, where $\delta \in (0,1)$, then we have, with probability $1-\delta$,
\[
	\sup_{A \in \A} |\mu_n(A) - \mu(A)| \leq 2 \sqrt{\frac{2p}{n} \big[ \log(12/\delta) + \log \S_\A(2n) 	\big]}.
\]
\end{theoreme}
\begin{proof}[Proof of Theorem~\ref{thm: relative deviations}]
It follows from~\citet[Theorem~1.11]{Lugosi2002} that
\begin{equation}
\label{eq: left deviation}
	\P \lp \sup_{A \in \A} \frac{\mu(A)-\mu_n(A)}{\sqrt{\mu(A)}} \geq t \rp \leq 4 \S_\A(2n) \e^{-nt^2/4}
\end{equation}
and
\begin{equation}
\label{eq: right deviation}
	\P \lp \sup_{A \in \A} \frac{\mu_n(A)-\mu(A)}{\sqrt{\mu_n(A)}} \geq t \rp \leq 4 \S_\A(2n) \e^{-nt^2/4},
\end{equation}
for every $t>0$. Therefore, it is convenient to decompose our probability of interest as follows
\begin{align}
	\P \lp \sup_{A \in \A} \frac{|\mu(A)-\mu_n(A)|}{\sqrt{p}} \geq t \rp
	&\leq \P \lp \sup_{A \in \A} \frac{\mu(A)-\mu_n(A)}{\sqrt{\mu(A)}} \geq t \rp + \P 								\lp \sup_{A \in \A} \frac{\mu_n(A)-\mu(A)}{\sqrt{p}} \geq t \rp \nonumber \\
	&\leq 4 \S_\A(2n) \e^{-nt^2/4} + \P \lp \mu_n(\mathbb{A}) \geq 2p \rp + \P \lp 							\sup_{A \in \A} \frac{\mu_n(A)-\mu(A)}{\sqrt{\mu_n(A)}} \geq t/\sqrt{2} \rp \nonumber \\
	&\leq 4 \S_\A(2n) \e^{-nt^2/4} + \P \lp \mu_n(\mathbb{A}) \geq 2p \rp + 4 								\S_\A(2n) \e^{-nt^2/8}, \label{eq: splitting mass}
\end{align}
where we made use of~\eqref{eq: left deviation} and ~\eqref{eq: right deviation}. By hypothesis, the second term of~\eqref{eq: splitting mass} is bounded by $\delta/3$. Indeed, by Bernstein's inequality~\eqref{eq: bernstein ineq},
\[
	\P \lp \mu_n(\mathbb{A}) \geq 2p \rp = \P \lp K -np \geq np \rp \leq \exp \lp - \frac{3}{8} np \rp,
\]
where $K \sim \Bin(n,p)$ as in Lemma~\ref{lem: cond trick}. Choosing $t=t(\delta)=2 \sqrt{\frac{2}{n} \big[ \log(12/\delta) + \log \S_\A(2n) 	\big]}$, we get that the first term and last term of~\eqref{eq: splitting mass} are also bounded by $\delta/3$. The result follows.
\end{proof}

%% file: proofSymBefore.tex
\begin{proof}[Proof~of~Theorem~\ref{thm: VC ineq rare sym before}]
	Let $a \in (0,1)$ and $t>0$. By symmetrization (Lemma~\ref{lem: symmetrization}), if $4(1-a)^2nt^2 \geq 2$,
	\[
	\P \lp \sup_{A \in \A} \left| \mu_n(A) - \mu(A) \right| \geq t \rp \leq 4 \P \lp \sup_{A \in \A} (\mu_n(A) - \mu_n'(A)) \geq at \rp,
	\]
	where $\mu_n'$ is the empirical measure associated to $X_1',\ldots,X_n'$, an independent ghost sample with common distribution $\mu$ and independent of $X_1,\ldots,X_n$. We easily verify that, since $\mu(A) \leq p$ for all $A \in \A$, the symmetrization remains true under the weaker hypothesis that 
	\begin{equation}
		\label{eq: weaker condition}
		(1-a)^2nt^2 \geq 2p.
	\end{equation}
	It suffices to adapt the bound on the variance $\Var[\1{\{X_i \in A\}}]$ in the proof of~Lemma~\ref{lem: symmetrization}.
	
	Let $\tilde{K} = \sum_{i=1}^n \1 {\{X_i \in \mathbb{A} \text{ or } X_i' \in \mathbb{A} \}} = \sum_{i=1}^n \1 {\{(X_i,X_i') \in \tilde{\mathbb{A}}\}}$, where $\tilde{\mathbb{A}} = (\mathbb{A} \times \R^d) \cup (\R^d \times \mathbb{A})$. Then, if $(Y_i,Y_i')_{i=1}^n$ is an independent sample, also independent of $(X_i,X_i')_{i=1}^n$, with common distribution given by the conditional distribution of $(X_1,X_1')$ when it lies in $\tilde{\mathbb{A}}$, we have by Lemma~\ref{lem: cond trick} 
	\begin{align*}
		&\P \lp \sup_{A \in \A} (\mu_n(A) - \mu_n'(A)) \geq at \rp \\
		&\qquad \qquad \qquad = \P \lp \sup_{A \in \A} \frac{1}{n} \sum_{i=1}^n \lp \1{\{X_i \in A \}} - \1{\{X_i' \in A \}} \rp \geq at \rp \\
		&\qquad \qquad \qquad = \E \lc \P \lp \sup_{A \in \A} \frac{1}{n} \sum_{i=1}^n \lp \1{\{X_i \in A \}} - \1{\{X_i' \in A \}} \rp \geq at \, \Big| \, \tilde{K} \rp \rc \\
		&\qquad \qquad \qquad = \E \lc \P \lp \frac{\tilde{K}}{n} \sup_{A \in \A} \frac{1}{\tilde{K}} \sum_{i=1}^{\tilde{K}} \lp \1{\{Y_i \in A \}} - \1{\{Y_i' \in A \}} \rp \geq at \, \Big| \, \tilde{K} \rp \rc \\
		&\qquad \qquad \qquad = \E \lc \P \lp \sup_{A \in \A} \sum_{i=1}^{\tilde{K}} \lp \1{\{Y_i \in A \}} - \1{\{Y_i' \in A \}} \rp \geq nat \, \Big| \, \tilde{K} \rp \rc. 
	\end{align*}
	
	Randomizing using \emph{Rademacher random variables}, we have
	\[
	\sum_{i=1}^{\tilde{K}} \lp \1{\{Y_i \in A \}} - \1{\{Y_i' \in A \}} \rp \stackrel{d}{=} \sum_{i=1}^{\tilde{K}} \sigma_i \lp \1{\{Y_i \in A \}} - \1{\{Y_i' \in A \}} \rp,
	\]
	where $\tilde{K}$ is fixed and $\sigma_1,\ldots,\sigma_{\tilde{K}}$ are i.i.d. Rademacher random variables independent of $(Y_1,Y_1'),\ldots,(Y_n,Y_n')$, i.e., for every $i \in \{1,\ldots, \tilde{K}\}$, $\P(\sigma_i = 1) = \P(\sigma_i = -1) = 1/2 $.
	Conditioning on $(Y_i,Y_i')_{i=1}^n$ and applying Hoeffding's inequality \cite[Theorem~1.2]{Lugosi2002} to each of the $\S_\A(2\tilde{K})$ possible vectors $(\1{\{Y_i \in A \}} - \1{\{Y_i' \in A \}})$ that can arise as $A$ ranges over $\A$, we obtain
	\begin{align*}
		\lefteqn{\P \lp \sup_{A \in \A} \sum_{i=1}^{\tilde{K}} \lp \1{\{Y_i \in A \}} - \1{\{Y_i' \in A \}} \rp \geq nat \, \Big| \, \tilde{K} \rp} \\
		&= \E_Y \lc  \P_\sigma \lp \sup_{A \in \A} \sum_{i=1}^{\tilde{K}} \sigma_i \lp \1{\{Y_i \in A \}} - \1{\{Y_i' \in A \}} \rp \geq nat \, \Big| \, \tilde{K},Y_1,Y_1',\ldots,Y_n,Y_n' \rp \rc \\
		&\leq \min\left\{ 1,\S_\A(2\tilde{K}) \exp \lp - \frac{a^2 (nt)^2}{2 \tilde{K}} \rp \right\},
	\end{align*}
	where the last inequality follows from the fact that the expectation $\E_Y$ is taken only with respect to the random variables $Y_1,Y_1',\ldots,Y_n,Y_n'$. We deduce
	\begin{equation}
		\label{eq: bound by min}
		\P \lp \sup_{A \in \A} \left| \mu_n(A) - \mu(A) \right| \geq t \rp \leq 4 \E \lc \min\left\{ 1,\S_\A(2\tilde{K}) \exp \lp - \frac{a^2 (nt)^2}{2 \tilde{K}} \rp \right\} \rc. 
	\end{equation}
	To be able to deal with this last quantity, we consider whether $\tilde{K}$ is less than $2np(1+s)$ or not, where $s>0$ is to be determined. By monotonicity, it gives
	\[
	\E \lc \min\left\{ 1,\S_\A(2\tilde{K}) \exp \lp - \frac{a^2 (nt)^2}{2 \tilde{K}} \rp \right\} \rc  \leq \P \big( \tilde{K} \geq 2np(1+s) \big) + \S_\A(4np(1+s)) \exp \lp - \frac{n a^2 t^2}{4p(1+s)} \rp. 
	\]
	Since $a \in (0,1)$ is not fixed, we may, given $n$ and $t$, choose the largest $a < 1$ such that the condition~\eqref{eq: weaker condition} is respected. It is given by $a = a(n,t) = 1-\sqrt{2p/(nt^2)}$. For such an $a$, we have
	\begin{equation}
		\label{eq: monotonicity}
		\E \lc \min\left\{ 1,\S_\A(2\tilde{K}) \exp \lp - \frac{a^2 (nt)^2}{2 \tilde{K}} \rp \right\} \rc
		\leq \P \big( \tilde{K} \geq 2np(1+s) \big) + \S_\A(4np(1+s)) \exp \lp - \frac{nt^2 - 2\sqrt{2np}t + 2p}{4p(1+s)} \rp. 
	\end{equation}
	Let $\delta \in (0,1)$. We choose $s = s(\delta)$ such that the first term of \eqref{eq: monotonicity} is bounded $\delta/8$. Then, we choose $t=t(\delta)$ such that the second term equals $\delta/8$. As a consequence of \eqref{eq: bound by min}, we will deduce that $\P \lp \sup_{A \in \A} \left| \mu_n(A) - \mu(A) \right| \geq t(\delta) \rp \leq \delta$, or, equivalently, $\sup_{A \in \A} \left| \mu_n(A) - \mu(A) \right| \leq t(\delta)$ with probability at least $1-\delta$. 
	
	By Bernstein's inequality \eqref{eq: bernstein ineq}, since $\tilde{K} \sim \Bin(n,\tilde{p})$ with $\tilde{p} = p(2-p)$, we have
	\begin{align*}
		\P \big(\tilde{K} \geq 2np(1+s) \big)  
		&= \P \big(\tilde{K} - n\tilde{p} \geq 2np(1+s) - n\tilde{p} \big) = \P \big(\tilde{K} - n\tilde{p} \geq np(2s+p) \big) \\
		&\leq \exp \lp - \frac{(np)^2(2s+p)^2}{2(n\tilde{p}(1-\tilde{p})+np(2s+p)/3)} \rp \\
		&\leq \exp \lp - \frac{np(2s+p)^2}{2((2-p)+(2s+p)/3)} \rp 
		\leq \exp \lp - \frac{nps^2}{1+s/3} \rp. 
	\end{align*}
	Equaling this expression to $\delta/8$, we obtain a quadratic equation in $s$ whose positive solution is given by
	\[
	s_+ = s(\delta) := \frac{\log(8/\delta)}{3np} + \sqrt{\frac{\log(8/\delta)}{np}}.
	\]
	Since we assumed $np \geq 2\log(8/\delta)$, we have $s(\delta) \leq \frac{1}{6} + \frac{1}{\sqrt{2}} \leq 1$. Hence, by taking $s=1$, the first term of \eqref{eq: monotonicity} is bounded by $\delta/8$. 
	
	For such a choice of $s$, the second term of \eqref{eq: monotonicity} becomes
	\[
	\S_\A(8np) \exp \lp - \frac{nt^2 - 2\sqrt{2np} t + 2p}{8p} \rp. 
	\]
	This last expression equals $\delta/8$ if
	\[
	t = t(\delta) := 2 \sqrt{\frac{2p}{n} \big[ \log(8/\delta) + \log \S_\A(8np) \big]} + \sqrt{\frac{2p}{n}}.   \qedhere
	\]  
\end{proof}